\documentclass[12pt]{article}
\usepackage{amssymb, amsthm, graphicx, verbatim}
\usepackage{epsf}

\newcommand{\bfz}{{\mathbb {Z}}}

\newcommand{\bfc}{{\mathbb {C}}}
\newcommand{\bfr}{{\mathbb {R}}}
\newcommand{\bfq}{{\mathbb {Q}}}
\newcommand{\dirac}{\partial \kern -.086in / }
\newcommand{\alphah}{\overline{\alpha}}

\begin{document}

\title{On Stein fillings of the 3-torus $T^3$}
\author{Andr\'{a}s I. Stipsicz}
\date{}
\maketitle

\newtheorem{fact}{Fact}[section]
\newtheorem{lemma}[fact]{Lemma}
\newtheorem{theorem}[fact]{Theorem}
\newtheorem{definition}[fact]{Definition}
\newtheorem{remark}[fact]{Remark}
\newtheorem{remarks}[fact]{Remarks}
\newtheorem{corollary}[fact]{Corollary}
\newtheorem{proposition}[fact]{Proposition}
\newtheorem{conjecture}[fact]{Conjecture}
\newtheorem{problem}[fact]{Problem}

\newenvironment{prooff}{\medskip \par \noindent {\it Proof}\ }{\hfill
$\square$ \medskip \par}
	\def\sqr#1#2{{\vcenter{\hrule height.#2pt
   		\hbox{\vrule width.#2pt height#1pt \kern#1pt
      		\vrule width.#2pt}\hrule height.#2pt}}}
	\def\square{\mathchoice\sqr67\sqr67\sqr{2.1}6\sqr{1.5}6}
\def\pf#1{\medskip \par \noindent {\it #1.}\ }
\def\endpf{\hfill $\square$ \medskip \par}
\def\demo#1{\medskip \par \noindent {\it #1.}\ }
\def\enddemo{\medskip \par}
\def\qed{~\hfill$\square$}

\def\pf#1{\medskip \par \noindent {\it #1.}\ }
\def\endpf{\hfill $\square$ \medskip \par}
\def\demo#1{\medskip \par \noindent {\it #1.}\ }
\def\enddemo{\medskip \par}
\def\qed{~\hfill$\square$}

\begin{abstract}
Topological properties of Stein fillings of contact 3-manifolds
diffeomorphic to the 3-torus 
$T^3$ are determined. We
show that for a Stein filling $S$ of $T^3$
the first
Betti number $b_1(S)$ is two, while $Q_S=\langle 0\rangle$. 
In the proof we also show that if $S$ is Stein and 
$\partial S$ is diffeomorphic to the Seifert fibered 3-manifold
$-\Sigma (2,3,11) $ then $b_1(S)=0$ and $Q_S=H$.  
Similar results are obtained for $S^1\times S^2$. Finally, we
describe Stein fillings of the Poincar\'e homology sphere 
$\pm \Sigma (2,3,5)$; in studying these
fillings we apply recent gauge theoretic results, and prove our
theorems by determining certain Seiberg-Witten invariants.
\end{abstract}

\section{Introduction}
\label{egy}
Suppose that $M$ is a closed, oriented 3-manifold and $\xi $ is a 2-plane
field on $M$. This $\xi$ is called a {\em contact structure\/}
if it is completely nonintegrable, i.e., for the 1-form (locally) defining
$\xi$ as $\ker \alpha$, the expression $\alpha \wedge d\alpha $ is 
nowhere 0. 
(For more about contact manifolds see \cite{Ae}.)
A 3-manifold $M$ in a K\"ahler surface $X$ inherits a natural contact 
structure provided that $M$ is {\em convex\/}, that is, 
there exists a vector field $v$ on $X$ transverse to $M$ such that
${\cal {L}}_v \omega =\omega$ for the K\"ahler form $\omega$
(${\cal {L}}$ stands for the Lie derivative).
In this case the complex lines in $TM$ form a 2-plane field $\xi$ satisfying 
the definition of a contact structure. For example, if $X$ admits a proper
biholomorphic embedding into $\bfc ^n$ for some $n$ (that is, $X$ is a
{\em Stein surface\/}) then for the distance function $f=\vert \vert \ .\ -p
\vert \vert ^2 
\colon X\to [0, \infty )$ (for $p\in \bfc ^n$ generic) the submanifolds
$f^{-1}(t)$ will be convex, hence contact away from the critical 
points. In fact, this property characterizes Stein surfaces:

\begin{theorem}[\cite{Gr}]\label{gra}
The (noncompact) complex surface $X$ is Stein if and only if there is a 
proper Morse function $f\colon X\to [0, \infty )$  such that away from 
the critical points the submanifolds $f^{-1}(t)\subset X$, with the 
plane fields induced by the complex tangent lines of $X$ in $TM$, are 
contact 3-manifolds. \qed    
\end{theorem}

Suppose that $t\in \bfr $ is a regular value of the above Morse function
$f\colon X\to [0,\infty )$. The manifold (with boundary) 
$S=f^{-1}[0,t]\subset X$ is called a {\em Stein domain\/}, and it can be
regarded as the compact version of Stein surfaces. (For more about
Stein surfaces and Stein domains see \cite{G, GS}.)

\begin{definition}
{\rm The contact manifold $(M, \xi )$ is {\em Stein fillable\/} if
there is a Stein domain $S$ such that $(M, \xi )$ is contactomorphic
to $\partial S$ (with the induced contact structure on it). In this
case $S$ is called a {\em Stein filling\/} of $(M, \xi )$.}
\end{definition}

\begin{remark}
{\rm
We always assume that $M$ is oriented and $\xi $ respects this orientation
through the requirement that $\alpha \wedge d\alpha >0$ for any 
1-form $\alpha$ defining $\xi $. Since $S$ has a natural orientation 
(as a complex surface), it induces an orientation on $\partial S$.
We require that the above contactomorphism is orientation 
preserving.}
\end{remark}

It is expected that the knowledge of all contact structures on $M$ will tell
us something about its geometry. To achieve this goal it seems reasonable 
to study all Stein fillings of a given 3-manifold. On the other hand, Stein 
domains can be regarded as analogues of minimal complex surfaces of 
general type in the category of manifolds with boundary. 
Therefore the study of Stein domains is interesting from the 4-dimensional
point of view as well. 
The {\em geography problem\/} for surfaces of general type asks  the 
possible values of $b_1$, $c_1^2$ and $c_2$ of such manifolds. Extending this 
problem we get:
\begin{problem}[The geography problem for Stein domains]
{\rm 
Fix a contact 3-manifold $(M,\xi )$  and describe characteristic 
numbers of Stein fillings of it.}
\end{problem}

In this paper we 
will address the problem of describing Stein domains with 
the 3-torus $T^3$, $S^1\times S^2$, $- \Sigma (2,3,11)$ and
$\pm \Sigma (2,3,5)$ as contact boundary. 
(For a possible definition of these Seifert fibered manifolds see
Figure~\ref{nuc}.)
The problem of Stein fillability
(and more generally, symplectic fillability) of contact 3-manifolds
has been extensively studied recently, see for example 
\cite{DG, E1, L1, L2, LM, Mc}. We only mention a prototype result here:

\begin{theorem}[\cite{E1}]\label{fill}
If $W$ is a Stein domain with $\partial W = S^3$ the 3-dimensional
sphere then $W$ is 
diffeomorphic to the 4-dimensional disk $D^4$. \qed
\end{theorem}

In the following we will prove a similar (but substantially weaker)
statement for the 3-torus $T^3$, $S^1\times S^2$ and for the Seifert
fibered 3-manifolds $- \Sigma (2,3,11)$ and $\pm \Sigma (2,3,5)$. Our
main result determines homological properties of Stein fillings of the
3-torus $T^3$. The intersection form of a 4-manifold
$X$ will be denoted by $Q_X$.

\begin{theorem}\label{torusz}
If $S$ is a Stein filling of $T^3$ then $b_1(S)=2$ and 
$Q_S=\langle 0 \rangle$; in particular, $\chi (S)=0$ and $\sigma (S)=0$.
Moreover, an appropriate Stein structure on $D^2\times T^2$ provides a 
Stein filling of $T^3$ with the above properties. 
\end{theorem}

\begin{remark}
{\rm In our proof we show that, in fact, $\pi _1(S)\cong \bfz \oplus
\bfz $ also holds. Using this isomorphism and Freedman's Classification
Theorem one can prove that the Stein filling $S$ is homeomorphic to 
$D^2\times T^2$. We hope to return to this question later, see 
Remark~\ref{homeo} and \cite{Sgok}.}
\end{remark}

\noindent
A similar argument will show
\begin{theorem}\label{s1s2}
If $S$ is a Stein filling of $S^1\times S^2$ then 
$\pi _1(S)\cong \bfz $ and $b_2(S)=0$.
\end{theorem}
\begin{remark}
{\rm
With a little  more care one can actually determine Stein fillings
of $\#  nS^1\times S^2$ up to diffeomorphism, and get that if 
$S$ is Stein and $\partial S\cong \#  nS^1\times S^2$ then $S$
is diffeomorphic to the $n$-fold boundary connected sum
of $S^1\times D^3$. For details see
\cite{Sgok}.}
\end{remark} 

The proof of Theorem~\ref{torusz} rests on the following result. 
\begin{theorem}\label{main1}
If $S$ is a Stein filling of $-\Sigma (2,3,11)$ for some contact 
structure on it, then $b_1(S)=0$ and $Q_S=H$. Moreover, there is a 
Stein domain $S$ with $b_1(S)=0$, $Q_S=H$ and $\partial S=-\Sigma (2,3,11)$.
\end{theorem}


(Here, as usual, $H$ denotes the hyperbolic plane 
$H=\left[{0\atop 1}\ {1\atop 0}\right]$.
Below $E_8$ stands for the symmetric bilinear form defined by the 
negative definite Cartan matrix of the exceptional Lie algebra $E_8$.
For the fixed orientation of the 3-manifolds $\Sigma (2,3,5)$ and
$\Sigma (2,3,11)$ see text 
following Theorem~\ref{roh}.) 
Similar arguments as applied in the
above theorem show the following --- these statements were already known 
\cite{AO2, L1, OO}.

\begin{theorem}\label{235}
If $S$ is a Stein filling of $\Sigma (2,3,5)$ then 
$b_1(S)=0$ and $Q_S=E_8$. The 3-manifold $-\Sigma (2,3,5)$ admits 
no Stein filling.
\end{theorem}
\begin{remark}
{\rm Notice the difference between {\em symplectic\/}
and {\em Stein\/} fillings. Blowing up is a symplectic operation which
ruins the Stein structure. Even if we disregard blow-ups, no finiteness
can be expected for symplectic fillings for $-\Sigma (2,3,11)$ or 
$T^3$ since by fiber summing the obvious fillings with elliptic surfaces 
infinitely many symplectic fillings with distinct $b_2^+$-invariant can
be constructed.}
\end{remark}

In computing the first Betti numbers of various Stein fillings we will
verify the following more general statement:
\begin{proposition}\label{begyes}
If $S$ is a Stein filling of a contact 3-manifold $(M, \xi )$ then 
the homomorphism $i_*\colon \pi _1 (M)\to \pi _1(S)$
induced by the inclusion $i\colon M\to S$ is a surjection.
Consequently $i_*\colon H_1(M;\bfz )\to H_1(S;\bfz )$ is onto, hence
$b_1(S)\leq b_1(M)$.
\end{proposition}

In proving Theorem~\ref{main1} we will use
recent results in gauge theory. The relevant theorems and constructions
will be summarized in Section~\ref{ketto}. 
Section~\ref{harom} deals with fillings of $-\Sigma (2,3,11)$ while
Section~\ref{haromfel} contains 
the proof of our main result Theorem~\ref{torusz}. In the final section 
we prove Theorem~\ref{235}.

\bigskip

\noindent 
{\bf {Acknowledgement}}: We thank Paolo Lisca for helpful comments
and Mustafa Korkmaz for drawing the figures.
During the course of this work the 
author was partially supported by OTKA T034885 and Sz\'echenyi
Professzori \"Oszt\"ond{\'\i}j.

\section{Gauge theoretic backgrounds}
\label{ketto}
We will frequently invoke the following celebrated result of Donaldson:
\begin{theorem}[\cite{DK}]\label{dona}
If $X$ is a smooth, closed 4-manifold with negative definite intersection 
form, then $Q_X$ is standard, that is, isomorphic to 
$\oplus _{1} ^{b_2(X)} \langle -1\rangle$. If $X$ is a  smooth, 
simply connected spin 4-manifold with $b_2^+(X)=1$ 
then $Q_X$ is isomorphic to $H$.
\qed
\end{theorem}
\begin{remark}\label{furu}
{\rm 
Using the monopole equations rather the instantons Donaldson 
originally used in his proof, Furuta \cite{F} extended Theorem~\ref{dona}
by showing that if a smooth spin 4-manifold $X$ has $Q_X=2kE_8\oplus 
lH$ then $l\geq 2\vert k \vert +1$.}
\end{remark} 
At one point we will appeal to the following famous result of Rohlin:
\begin{theorem}[\cite{Ro}] \label{roh}
If $X$ is a smooth spin 4-manifold then the signature $\sigma (X)$
of $X$ is divisible by 16.
\qed
\end{theorem}

The 3-manifolds $\Sigma (2,3,5)$ and $\Sigma (2,3,11)$ are defined 
as oriented boundaries of the complex manifolds $M_c(2,3,5)=\{ 
(x,y,z) \in \bfc ^3 \mid x^2+y^3+z^5=\varepsilon ,\ \vert x\vert ^2 +
\vert y\vert ^2 +\vert z\vert ^2\leq 1 \}$ and
$M_c(2,3,11)=\{ 
(x,y,z) \in \bfc ^3 \mid x^2+y^3+z^{11}=\varepsilon ,\ \vert x\vert ^2 +
\vert y\vert ^2 +\vert z\vert ^2\leq 1 \}$ (with $\vert \varepsilon \vert 
\ll 1 $), i.e., as the boundaries of the corresponding (compactified)
Milnor fibers.
\begin{figure}
\centerline{		
\epsfbox{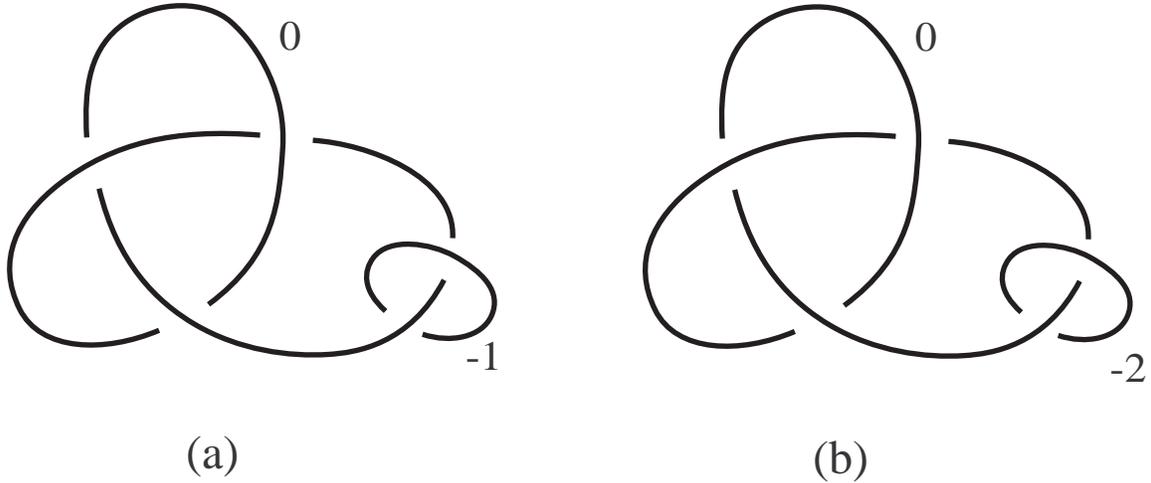}}
\caption{Kirby diagrams for (a) $N_1$ and (b) $N_2$}		
\label{nuc}
\end{figure}
Alternatively, $-\Sigma (2,3,5) $ and $-\Sigma (2,3,11)$ are the 
oriented boundaries of the 
the nuclei $N_1\subset E(1)$ and $N_2\subset E(2)$, 
where $N_1$ and $N_2$ are given 
by the Kirby diagrams of Figure~\ref{nuc}.

In our subsequent arguments we will apply product formulae for
Seiberg-Witten invariants when we cut 4-manifolds along 
$\Sigma (2,3,5)$ or $\Sigma (2,3,11)$. For the sake of completeness we sketch
the definition of Seiberg-Witten invariants and analyze the 3-dimenisonal
equations for $\Sigma (2,3,5)$ and $\Sigma (2,3,11)$ more carefully.

Suppose that $X$ is a smooth, closed, oriented 4-manifold with no 2-torsion
in $H_1(X; \bfz )$. A characteristic element $K$ in 
$H^2(X; \bfz )$ uniquely determines a spin$^c$ structure on $X$, and once
a connection $A$ on the line bundle $L$ with $c_1(L)=K$ is fixed, a metric on
$X$ gives rise to a twisted Dirac operator $\dirac _A$. 
The space of connections on the complex line bundle $L$
will be denoted by ${\cal {A}}_L$, the curvature of
$A\in {\cal {A}}_L$ is $F_A$, and $F_A^+$ stands for its self-dual part.
Let $W^+$ denote the positive spinors corresponding to the spin$^c$ structure
defined by $K$. Then the Seiberg-Witten equations for a pair
$(A, \psi ) \in {\cal {A}}_L\times \Gamma (W^+)$  read as
$$\dirac _A \psi =0\ \ \ \  {\rm {and}}\ \ \ \ \   
F_A^+=iq (\psi ). \leqno(\ast ) $$
(Here $q \colon \Gamma (W^+)\to \Gamma (\Lambda ^+)$ is a certain 
quadratic map.) 
If $K^2-(3\sigma (X)+\chi (X))=0$, the solution 
space ${\cal {M}}_K$ 
of the equations (mod symmetries of the equations) is a compact 
0-dimensional manifold. An orientation of 
$H_+^2(X; \bfr )\otimes H^1(X;\bfz )$ fixes an orientation on this
solution space, and provided $b_2^+(X)>1$, the algebraic sum of the points 
of ${\cal {M}}_K$ turns out to be a smooth invariant of $X$, denoted
by $SW_X(K)\in \bfz $. Similar, but somewhat more complicated 
procedure provides $SW_X(K)\in \bfz $ for $K$ with 
$K^2-(3\sigma (X)+2\chi (X))>0$. (If 
$K^2-(3\sigma *X)+2\chi (X))<0$ or $K$ is not characteristic, then 
$SW_X(K)=0$ by definition.)
The class $K$
is called a {\em basic class\/} of $X$ if $SW_X(K)\neq
0$. It can be shown that 
$SW_X(K)= (-1)^{{\sigma (X)+\chi (X)}\over 4} SW_X(-K)$, therefore
$K$ and $-K$ are basic classes at the same time.
We say that $X$ is of {\em simple type\/} if for a basic class $K$
of $X$ the equation $K^2=3\sigma (X)+2\chi (X)$ holds.
Following ideas of Fintushel and Stern \cite{FS} we can associate 
a formal series to any 4-manifold $X$ with $b^+_2(X)>1$: 
If $\lbrace \pm K_1,...,\pm K_n\rbrace $ are the nonzero
basic classes of $X$, then take
$${\cal {SW}}(X)=SW_X(0)+\sum _{i=1}^n SW_X(K_i)\exp (K_i)+ SW_X(-K_i)
\exp (-K_i).$$
It has been proved \cite{W} that if $X$ is a minimal surface of general type
then it has two basic classes $\pm c_1(X)$, and $SW_X (\pm c_1(X))=\pm 1$.
Therefore in that case ${\cal {SW}}(X)=\exp (c_1(X))\pm \exp (-c_1(X)).$
(For other complex surfaces we only know that $\pm c_1(X)$ are basic classes.)
It is also known that the K3-surface $Y$ has a single basic class
which is $c_1(Y)=0$, hence ${\cal {SW}}(Y)=1$.
For a more thorough study of Seiberg-Witten theory, see \cite{GS, Mo, W}.

In a similar vein the 3-dimensional analogue of Seiberg-Witten
equations can be defined. 
For a 3-manifold $M$ the spin$^c$ structures are parametrized by
$H^2(M;\bfz )$ and if $W\to M$ denotes the spinor bundle then 
the Seiberg-Witten equations  for  $(A,\psi )\in 
{\cal {A}}_{\det W}\times \Gamma (W)$ are
$$\dirac _A \psi =0\ \ \ {\rm {and}}\ \ \ \ast F_A=iq (\psi ).\leqno(\ast
 \ast )$$
(As usual, $F_A$ denotes the curvature 2-form of the connection
$A\in {\cal {A}}_{\det W}$, 
and $\ast $ stands for the Hodge $\ast$-operator given
by a metric on $M$. Now $q$ maps from 
$\Gamma (W)$ to $\Gamma (\Lambda ^1 M)$.) These equations have been solved for
$\Sigma (2,3,11)$ in \cite{MOY}. Notice that since $\Sigma (2,3,11)$
is an integral homology sphere, it admits a unique spin$^c$ structure.
After substituting the Levi-Civita connection with a suitable connection
in the definition of $\dirac _A$, 
in \cite{MOY} it was shown that
$(\ast \ast)$ admits 3 solutions (up to gauge equivalence): 
one of them is
the trivial solution $\theta$, which is the trivial connection with vanishing
spinor field; the other two will be denoted by   $\alpha$ and
$\alphah$. 

\begin{remarks}
{\rm 
\begin{itemize}
\item
Such a perturbation of the Seiberg-Witten equations over a three-manifold 
can be naturally extended to give a perturbation of the Seiberg-Witten 
moduli space over 4-manifolds containing long necks.
It is proved  \cite{OSz} that this perturbation over a smooth
closed 4-manifold with $b_2 ^+ >1$ gives a compact moduli space 
which is smoothly cobordant to the unperturbed Seiberg-Witten 
moduli space. This implies that such a perturbation can be used to compute
the Seiberg-Witten invariants.

\item
Because of the presence of a positive scalar curvature metric, one can 
easily show that the Seiberg-Witten equations on $\Sigma (2,3,5)$ admit 
a unique solution $\theta$, which is the trivial connection with vanishing 
spinor field.
\end{itemize}
}
\end{remarks}

By finding relations between the $L^2$ moduli spaces of 
Seiberg-Witten solutions over a 4-manifold $X$ with boundary diffeomorphic to 
$\Sigma (2,3,11)$, in \cite{SSz} relative invariants, relative basic
classes and the (formal) series
${\cal {SW}}(X)$ has been defined for a smooth 4-manifold $X$ with
boundary diffeomorphic to $\pm \Sigma (2,3,11)$.
The relation between absolute
and relative invariants is given by 

\begin{theorem}[\cite{SSz}] \label{szorz} 
If the closed 4-manifold $Z$ decomposes as $Z=X\cup _{\Sigma (2,3,11)}Y$
with $b_2^+(X), b_2^+(Y)>0$
then 
${\cal {SW}}(Z)={\cal {SW}}(X)\cdot {\cal {SW}}(Y)$, that is,
the product of the relative invariants equals the absolute invariant of the
closed 4-manifold $Z$. \qed
\end{theorem}

There are three more important ingredients of the proofs we will describe in 
the following sections. The first theorem (due to Lisca and Mati\'c) 
provides a K\"ahler embedding of a Stein domain into a minimal surface
of general type, more precisely

\begin{theorem}[\cite{LM}]\label{lm}
For a Stein domain $S$ there exists a minimal surface $X$ of general type
and a K\"ahler embedding $f\colon S \to X$.  Moreover, we can assume
that $X-f(S)$ is not spin and
$b_2^+(X-f(S))>1$. \qed
\end{theorem}

The next theorem is a special case of a result of Ozsv\'ath
and Szab\'o which describes restrictions on the embedding of certain 
circle bundles over surfaces. Suppose that $M_{e,1}$ is a circle bundle
over the torus with Euler number $e$.
\begin{theorem}[\cite{OSz2}] \label{fish}
If the minimal surface $X$ of general type decomposes as
$X=X_1\cup _{M_{e,1}} X_2$ along the 3-manifold 
$M_{e,1}$ with $\vert e\vert \geq 1$ then 
either $b_2^+(X_1)=0$ or $b_2^+(X_2)=0$. \qed
\end{theorem}
Combinig Theorems~\ref{lm} and \ref{fish}  we  get (see also \cite{AO2}):
\begin{corollary}\label{b+0}
If $S$ is a Stein domain with $\partial S=M_{e,1}$ and $\vert e \vert
\geq 1$ then $b_2^+(S)=0$.\qed
\end{corollary}

Finally we invoke a result of Morgan and Szab\'o which characterizes 
homotopy K3-surfaces:
\begin{theorem}[\cite{MSz}]\label{msz}
Suppose that $X$ is a simply connected spin 4-manifold of simple type
with $Q_X=2kE_8\oplus lH$ and $SW_X(0)=\pm 1$. Then $Q_X=2E_8\oplus 3H$.
\qed
\end{theorem}

(Notice that since $X$ is of simple type and 0 is a basic class, it follows
that $l=4k-1$. The theorem of Morgan and Szab\'o proves that $SW_X(0)$ is even
once $k>1$.)

\section{Fillings of $-\Sigma (2,3,11)$}
\label{harom}

We begin our study of Betti numbers of Stein fillings by proving an estimate
on $b_1$.

\begin{prooff}{\em of Proposition~\ref{begyes}.}
It is well-known that a Stein domain $S$ can be built up using
0-, 1- and 2-handles only; for such manifolds the 
surjectivity of $\pi _1 (\partial S)\to \pi _1 (S)$ is obvious.
This surjection now trivially implies that 
$H_1(\partial S; \bfz )\to H_1(S; \bfz )$ is also a surjection, hence
$b_1 (S)\leq b_1(\partial S)$.
\end{prooff}

Since $\pm \Sigma (2,3,5)$ and $- \Sigma (2,3,11)$ are 
integral homology spheres, the above theorem shows that 
Stein fillings of these Seifert fibered 3-manifolds have trivial 
first homology, hence vanishing first Betti number.

Next we consider intersection forms of fillings of $- \Sigma (2,3,11)$.

\begin{proposition}\label{nonegdef}
If $X$ is a smooth 4-manifold with $\partial X=-\Sigma (2,3,11)$ then 
$b_2^+(X)>0$.
\end{proposition}
\begin{proof}
It is a standard fact that the K3-surface $Y$ contains three disjoint
copies of the nucleus $N_2$ (recall that $\partial N_2=-\Sigma
(2,3,11)$); and the intersection form of the manifold $Y-3{\mbox { int
}}N_2$ is negative definite and nonstandard.  So if $X$ is negative
definite with $\partial X = -\Sigma (2,3,11)$ then the 4-manifold
$(Y-3 N_2)\cup 3X$ we get by replacing the nuclei in $Y$ by $X$ is a
smooth manifold with nonstandard negative definite intersection
form. The existence of such a manifold, however, contradicts
Donaldson's Theorem~\ref{dona}, showing that $X$ is not negative
definite.
\end{proof}

\begin{prooff}{\em of Theorem~\ref{main1}.}
Let $S$ be a Stein filling of $-\Sigma (2,3,11)$ and 
consider the K\"ahler embedding $S\to X$ where $X$ is a minimal surface
of general type --- the existence of such an embedding is guaranteed by 
Theorem~\ref{lm}. Recall that we can assume that
$X-S$ is nonspin with $b_2^+(X-S)>1$. The product formula of 
Theorem~\ref{szorz} shows that ${\cal {SW}}(S)=\pm 1$: We use the
fact that $X$, as a minimal surface of general type has only two basic 
classes $\pm c_1(X)$; therefore ${\cal {SW}}(X)$ is nondivisible, but since
$X-S$ is nonspin, 0 is not characteristic and so ${\cal {SW}}(X-S)\neq 1$.
Notice that this computation shows that 
$c_1(S)=0$ is the unique basic class for $S$, in particular, $S$ is 
spin.  (The product formula of Theorem~\ref{szorz} applies since 
$b_2^+(S)>0$ by Proposition~\ref{nonegdef} 
and $b_2^+(X-S)>0$ by Theorem~\ref{lm}.)
Now consider $Z=(Y-N_2)\cup S$ --- as before, $Y$ stands for
the K3-surface. The product formula shows that ${\cal {SW}}(Z)=\pm 1$, and 
easy handle calculus verifies that $Z$ is simply connected:
$Y-N_2$ is simply connected and we can build $Z$ on the top of
it by adding only 2-, 3- and 4-handles, since $S$ is Stein.
In conclusion, for the simply connected spin manifold $Z$ we have 
that $SW_Z(0)=\pm 1$; applying Theorem~\ref{msz} of Morgan and Szab\'o
this fact implies that $Q_Z=2E_8\oplus 3H$ 
and since $Q_{Y-N_2}=2E_8\oplus 2H$, we get 
that $Q_S\cong Q_{N_2}=H$.
\end{prooff}
\begin{remarks}
{\rm
\begin{enumerate}
\item Figure~\ref{stein} demonstrates that, in fact, 
$N_2$ carries a Stein structure. (For handle calculus of
Stein domains, see \cite{G, GS}.) This provides a Stein filling of
$-\Sigma (2,3,11)$ as stated.
\begin{figure}
\centerline{		
\epsfbox{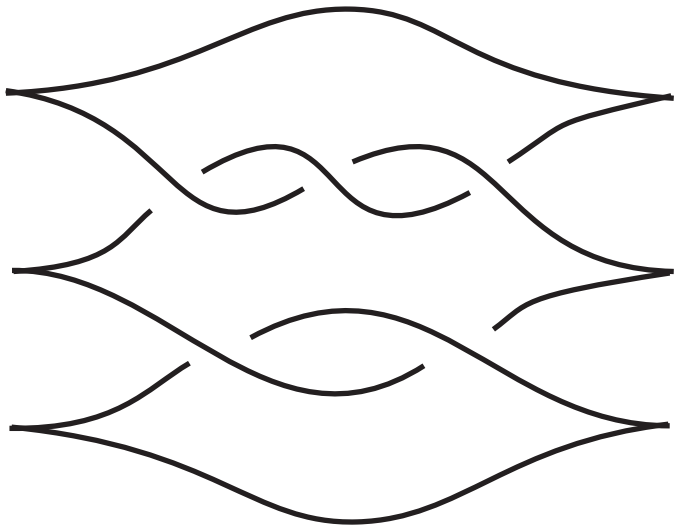}}
\caption{Stein structure on $N_2$}		
\label{stein}
\end{figure}
\item It is known \cite{Ho} that $-\Sigma (2,3,11)$ carries a unique
Stein fillable (in fact, a unique tight) contact structure.
\item It seems natural to conjecture that the Stein filling $S$ of
$-\Sigma (2,3,11)$ is diffeomorphic to $N_2$, although the techniques
applied in the above proof seem to be weak to verify such a conjecture.
\end{enumerate}
}
\end{remarks}

\section{Stein fillings of $T^3$}
\label{haromfel}
Using a famous result of Eliashberg \cite{E2} together with 
the classification result due to Giroux and Kanda
(given below),
now we can prove our result regarding Stein fillings of the 
3-torus $T^3$. Our proof of Theorem~\ref{torusz} will heavily rely on
Theorem~\ref{main1}.

\begin{figure}
\centerline{		
\epsfbox{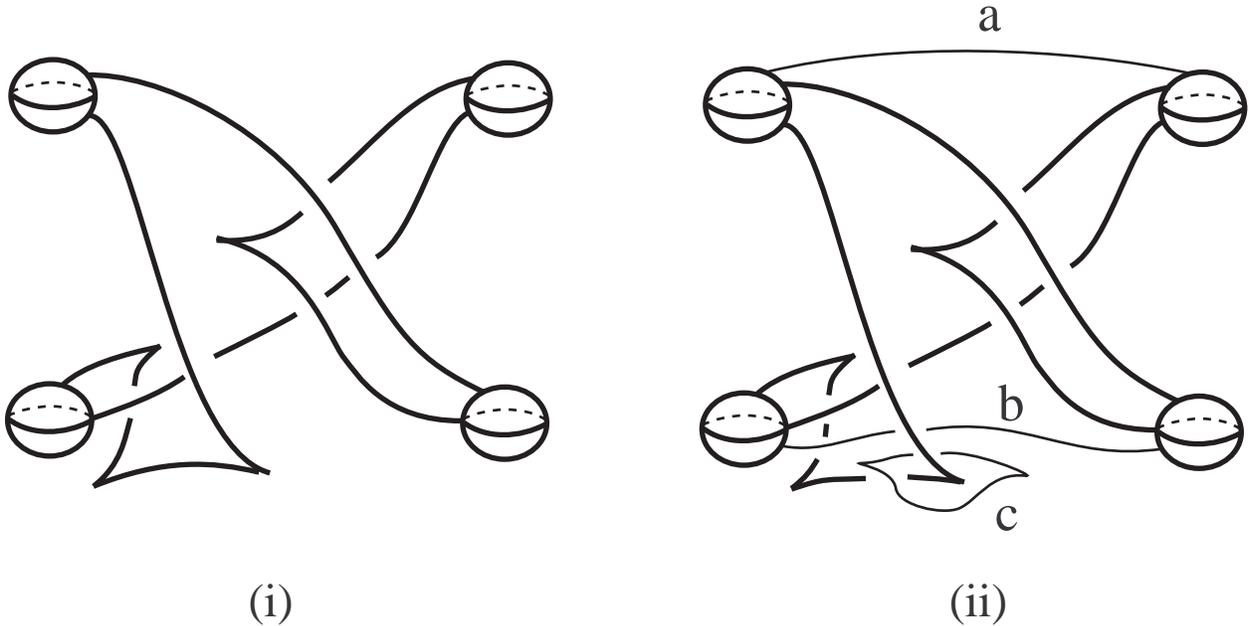}}
\caption{Stein domain with boundary $T^3$}
\label{t3}
\end{figure}

\begin{figure}
\centerline{		
\epsfbox{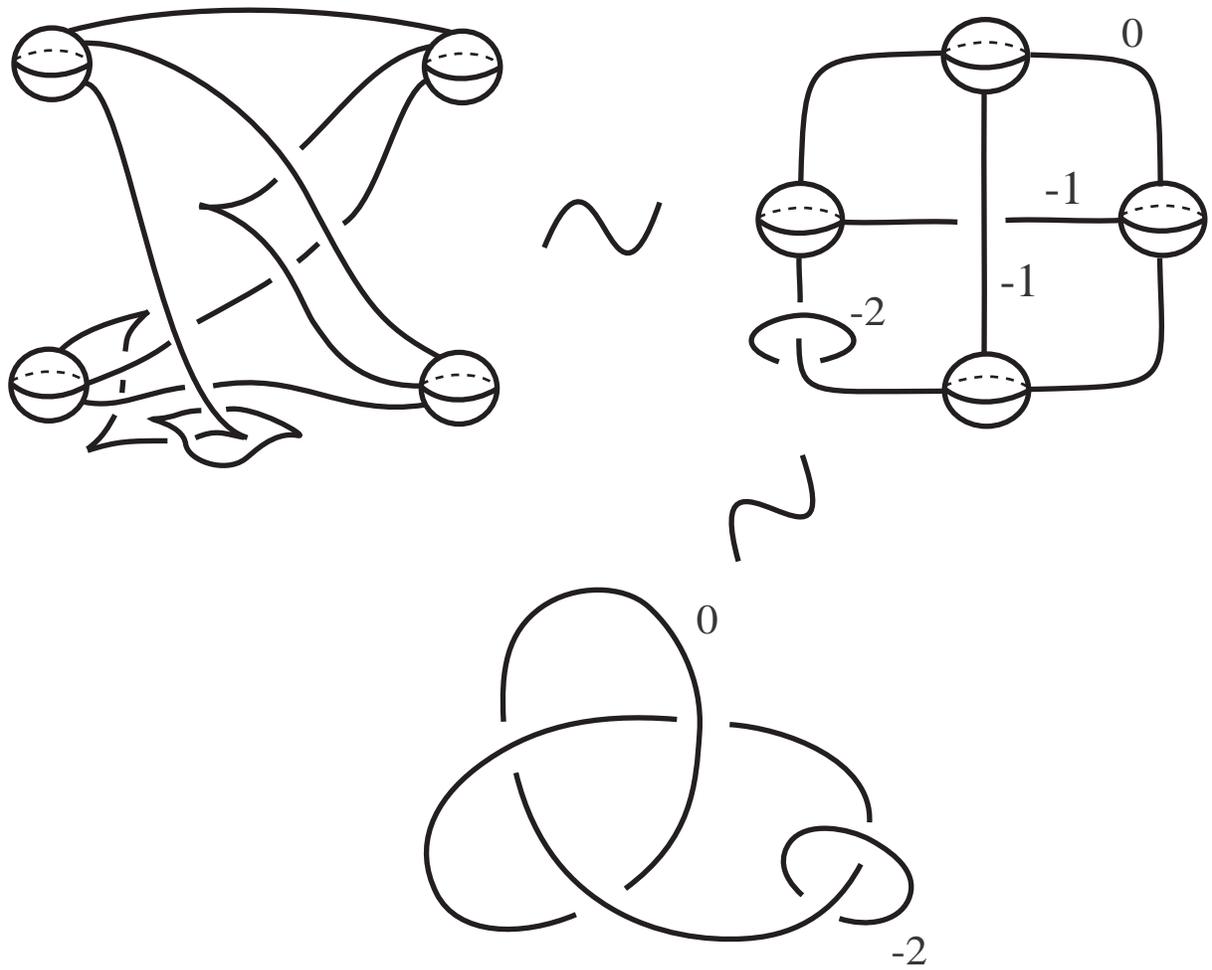}}
\caption{Proof of Lemma~\ref{fish-sig}(b) with Kirby diagrams}
\label{kirby}
\end{figure}

\begin{theorem}[\cite{E2, Gi, Ka}]\label{egk}
The contact structures $\xi _n =\ker (\cos (2\pi nt)dx-\sin (2\pi n t)dy)$
on $T^3$ $($in coordinates $(x,y,t)$ on $T^3)$
are all noncontactomorphic and comprise a complete list of tight 
contact structures on the 3-torus $T^3$. If $(T^3, \xi _n)$ admits a Stein 
filling then $n=1$. 
\qed
\end{theorem}
The contact structure $\xi _1$ can be given as the 
boundary of the Stein domain given by Figure~\ref{t3}(i).
(For the relation of Kirby diagrams and Stein structures
and for the definition of the Thurston-Bennequin invariant tb of a Legendrian 
knot, see \cite{G, GS}.) Before turning to the proof of Theorem~\ref{torusz}
we need a lemma. (The proof of it is standard Kirby calculus and is
left for the reader, see Figure~\ref{kirby}.)
\begin{lemma}\label{fish-sig}
(a) Gluing a 2-handle along the fine curve a (or b) of Figure~\ref{t3}(ii) 
with framing ${\mbox {tb}}(a)-1=-1$ results a handlebody 
with boundary $M_{1,1}$ of Theorem~\ref{fish}.

\noindent
(b) Gluing 2-handles along the three fine curves
a,b,c of Figure~\ref{t3}(ii) (with framings $-1, -1$ and $-2$, resp.) 
results a handlebody with boundary 
$-\Sigma (2,3,11)$.\qed
\end{lemma}

\begin{prooff}{\em of Theorem~\ref{torusz}.}
Fix a contactomorphism between $\partial S$ and the boundary of the
Stein domain of Figure~\ref{t3}(i).  For determining $\chi (S)$
consider $W=S\ \cup{\mbox{ three 2-handles}}$ attached along $a,b,c$
with framing ${\mbox{tb}}-1$. The result is a Stein filling of
$-\Sigma (2,3,11)$, and according to Theorem~\ref{main1} it has Euler
characteristic 3.  (Since we glue the 2-handles along Legendrian
curves with framing tb$-1$, the existence of a Stein structure on $W$
follows from by now standard arguments discussed in \cite{G,GS}.)
Removing the three 2-handles from $W$ 
we arrive to the conclusion $\chi (S)=0$.
This fact implies that $b_1(S)\geq 1$, since $b_1(S)=0$ and $\chi
(S)=0$ would imply $b_2(S)=-1$. Therefore $S$ admits unramified covers
of any degree. Notice that since $b_1(S)\leq 3$ by
Theorem~\ref{begyes} (and so $b_2(S)\leq 2$), we get that $\vert
\sigma (S)\vert \leq 2$. Consider a 3-fold unramified cover of $S$ ---
the result is a Stein filling ${\overline {S}}$ of the 3-fold cover of
$\partial S$, which is $T^3$ again. Therefore all the above said ---
in particular $\vert \sigma ({\overline {S}})\vert \leq 2$ --- holds
for ${\overline {S}}$. Since $\sigma ({\overline {S}})= 3\sigma (S)$,
we conclude that $\sigma (S)=0$.

Finally we show that $b_1(S)=2$ for a Stein filling $S$ of $T^3$.
According to Theorem~\ref{begyes} the map $H_1(T^3; \bfz )\to
H_1(S;\bfz )$ is onto --- we claim that the (image of the) circles $a$
and $b\subset T^3$ of Figure~\ref{t3}(ii) remain essential in $H_1(S;
\bfq )$ while $c$ becomes 0.  Since $\xi _n$ can be given as a
pull-back of $\xi _1$ under unramified cover along the $t$ coordinate,
if $i_*(c) \neq 0$ in $\pi _1(S)=H_1(S; \bfz )$ then a corresponding
$n$-fold cover of $S$ provides a Stein filling of $T^3$ equipped with
$\xi _n$. This contradicts Theorem~\ref{egk} once $n>1$, hence $c=0$
in $H_1(S;\bfz )$.  If $a=0$ in $H_1(S; \bfq )$ then attaching a
2-handle along $a$ with ${\mbox {tb}}(a)-1$ we get a Stein domain
${\tilde {S}}$ with $\partial {\tilde {S}}=M_{1,1}$ of
Lemma~\ref{fish} and $b_2^+({\tilde {S}})>0$, since the surface in $S$
with boundary $a$ together with the core of the handle and the dual
torus of $a$ in $T^3$ give a hyperbolic pair in ${\tilde {S}}$.  This
fact contradicts Corollary~\ref{b+0}, therefore $a\neq 0$ in $H_1(S;
\bfq)$.  The role of $b$ is analogous, hence $b\neq 0$ in $H_1(S;
\bfq)$.  It follows now that $S$ admits a CW decomposition with two
1-cells, and so the number of 2-cells is one in this decomposition
(since $\chi (S)=0$). Since $\pi _1(S)$ is Abelian (being a factor of
$\pi _1(T^3)\cong \bfz ^3$), we get that the attaching circle of this
2-cell is homologically trivial, therefore $H_1(S;\bfz ) \cong \pi
_1(S)\cong \bfz \oplus \bfz $. In particular, $b_1(S)=2$.  This
implies $b_2(S)=1$, and so the intersection form $Q_S$ can be easily
identified with $Q_S=\langle 0\rangle $. The proof is now complete.
\end{prooff}

\begin{remark}\label{homeo}
{\rm 
Notice that the above proof, in fact, showed that $\pi _1(S)\cong 
H_1(S;\bfz )$ is isomorphic to $\bfz \oplus \bfz$;
moreover we proved that a Stein filling of 
$T^3$ embeds into a Stein filling of $-\Sigma (2,3,11)$. 
Since for this latter 3-manifold all Stein fillings are spin, 
we conclude that a Stein filling of $T^3$ is spin. 
Since $\sigma (S)=0$, it follows that the induced
spin structure on $\partial S$ is diffeomorphic to the the one 
$\partial (Y-\nu F )$ inherits from $Y-\nu F$; here 
$Y$ is the K3-surface and $F$ is a regular fiber in an elliptic
fibration on $Y$. Therefore $Z=S\cup _{T^3}(Y-\nu F)$ is a 
(simply connected) spin 4-manifold, and using the gauge theoretic
results discussed in \cite{T} 
one can verify that $SW_Z(0)$ is $\pm 1$. Consequently 
Theorem~\ref{msz} implies that $Z$ is homeomorphic to $Y$.
Motivated by this homeomorphism one can conjecture that the
Stein filling $S$ is diffeomorphic to $D^2\times T^2$ ---
the Stein filling shown by Figure~\ref{t3}(i).
In fact, the extension of Freedman's Classification Theorem
for 4-manifolds with boundary and nontrivial fundamental group
shows that a Stein filling $S$ of $T^3$ is {\em homeomorphic\/}
to $D^2\times T^2$, see \cite{Sgok}.
}
\end{remark}
We close this section with the proof of Theorem~\ref{s1s2}.

\begin{prooff}{\em of Theorem~\ref{s1s2}.}
Since $S^1\times S^2$ admits a unique tight contact structure, 
adding a 1-handle and a 0-framed 2-handle (as shown by Figure~\ref{t3}(i))
to $\partial S$ results a Stein filling $W$ of $T^3$. Since $W=S\cup 
{\mbox {1-handle}}\cup {\mbox {2-handle}}$, we obviously have 
$\chi (S)=\chi (W)=0$. This again shows that $b_1(S)\geq 1$, which implies
$\pi _1(S)\cong \bfz $, since --- according to Proposition~\ref{begyes} ---
the fundamental group $\pi _1(S)$ is the factor of $\pi _1(S^1\times S^2)
\cong \bfz $. Now $b_2(S)=0$ trivially follows. 
\end{prooff}

\begin{remark}
{\rm
Using Theorem~\ref{fill} of Eliashberg one can also show that 
a Stein filling of $S^1\times S^2$ is diffeomorphic to $S^1\times D^3$,
see \cite{Sgok}. Similar argument determines the Stein filling of
$\# \ nS^1\times S^2$ for all $n\geq 1$ up to diffeomorphism.
}
\end{remark}

\section{Appendix: Stein fillings of $\pm \Sigma (2,3,5)$}
\label{negy}

Finally we turn to the proof of Theorem~\ref{235}.
The theorem was already proved by various authors 
(see Remark~\ref{4.5}); we include it here because the proof given 
below is very similar in spirit to the proof of Theorem~\ref{main1}. 
We begin our proof with the following ``folk'' theorem.
(For a complete proof see \cite{L1, OO}.)
\begin{lemma}\label{folk}
If $S$ is a Stein filling of a 3-manifold which admits positive 
scalar curvature then $S$ is negative definite.
\end{lemma}
\begin{prooff}{\em (sketch).}
Consider the K\"ahler embedding $S\to X$ where $X$ is a minimal surface
of general type with $b_2^+(X-S)>1$. Since a 3-manifold with positive scalar
curvature metric cannot divide a 4-manifold with nonzero Seiberg-Witten  
invariants into two pieces both with $b_2^+>0$, the lemma follows.
\end{prooff}
Notice that $\pm \Sigma (2,3,5)$ (as the quotient of $S^3$) admits a 
metric with positive scalar curvature.

\begin{proposition}
If $S$ is a Stein filling of $\Sigma (2,3,5)$ (for some contact structure)
then $S$ is negative definite and spin.
\end{proposition}
\begin{proof}
The fact that $b_2^+(S)=0$ follows from Lemma~\ref{folk}.
Recall that the Seiberg-Witten equations admit a unique
(up to gauge equivalence) solution on 
$\Sigma (2,3,5)$. Now consider an embedding $S\to X$ where 
$X$ is a minimal surface of general type.
Grafting solutions for the spin$^c$ structures $\pm c_1(S)$ and
$c_1(X-S)$ together we get 4 basic classes unless $c_1(S)=0$
or $c_1(X-S)=0$. Since $X$ has exactly two basic classes and
$c_1(X-S)$ is nonspin (therefore $c_1(X-S)\neq 0$) we get that
$c_1(S)=0$, consequently $S$ is spin.
\end{proof}

\begin{prooff}{\em of Theorem~\ref{235}.}
It can be easily verified that the negative definite
$E_8$-plumbing $E$ (see
Figure~\ref{e8}) embeds in the K3-surface $Y$.
\begin{figure}
\centerline{		
\epsfbox{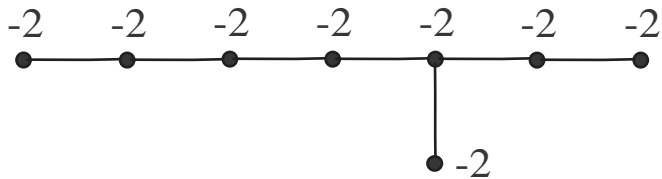}}
\caption{The negative definite $E_8$-plumbing}
\label{e8}
\end{figure}
For a Stein filling $S$ consider $Z=S\cup (Y-E)$. Since $S$ is spin
(and $\partial S =\Sigma (2,3,5)$ admits a unique spin structure) it is
spin and $\pi _1(S\cup (Y-E))=1$. 
Since $S$ is negative definite, we have that $Q_Z=(k+1)E_8\oplus 3H$.
Furuta's Theorem (see Remark~\ref{furu}) shows that $k\leq 1$; since
Rohlin's famous Theorem~\ref{roh} excludes $k=0$, we conclude that 
$Q_S=E_8$, proving Theorem~\ref{235} for $\Sigma (2,3,5)$. The part of
Theorem~\ref{235} about $-\Sigma (2,3,5)$ follows from the fact that 
$-\Sigma (2,3,5)$ does not bound negative definite 4-manifold at all:
If $\partial X=-\Sigma (2,3,5)$ and $b_2^+(X)=0$ then 
$X\cup _{-\Sigma (2,3,5)}E$ violates Donaldson's Theorem~\ref{dona}.
\end{prooff}

\begin{remarks}\label{4.5}
{\rm 
\begin{enumerate}
\item A similar proof of the above theorem was already found by Ohta and
Ono \cite{OO}. 
\item In the above proof one can also argue that a negative
definite spin 4-manifold with boundary diffeomorphic to 
$\Sigma (2,3,5)$ has intersection form isomorphic to $E_8$ along the
lines developed by Fr\o yshov, see \cite{Fro}.
\item According to \cite{Ho}, the 3-manifold $\Sigma (2,3,5)$
admits a unique fillable (in fact, a unique tight) contact structure.
Notice that in our proof we made no assumption on the contact structure
on $\Sigma (2,3,5)$.
\item The $E_8$-plumbing $E$ supports a Stein structure
providing a filling for $\Sigma (2,3,5)$.
It seems natural to expect that any Stein filling of 
$\Sigma (2,3,5)$ is diffeomorphic to $E$.
\item Using  Theorem~\ref{begyes} 
and analyzing possible quotients of $\pi _1(\Sigma (2,3,5))$ we get that
for a Stein filling $S$ of $\Sigma (2,3,5)$ the group $\pi _1(S)$ is trivial.
Therefore Freedman's Theorem implies that $S$ is, in fact,  homeomorphic
to the $E_8$-plumbing $E$. For the details of the above argument see 
\cite{Sgok}.
\item For $-\Sigma (2,3,5)$ Lisca proved that it
does not admit a symplectic semi-filling, while
Etnyre and Honda showed \cite{EH} that it supports no tight contact 
structure at all.
\end{enumerate}
}
\end{remarks}

\noindent
Department of Analysis, ELTE TTK,
1055. Kecskem\'eti u. 10-12., Budapest, Hungary\\
stipsicz@cs.elte.hu 


\begin{thebibliography}{AA}

\bibitem{Ae}
B. Aebisher et al., 
{\em Symplectic geometry},
Progress Math. {\bf124}, Birkh\"auser, Boston, MA, 1994.

\bibitem{AO1}
S. Akbulut and B. Ozbagci, 
{\em Lefschetz fibrations on compact Stein surfaces}, Geometry and Topology,
{\bf5} (2001), 319--334.


\bibitem{AO2}
S. Akbulut and B. Ozbagci, 
{\em On the topology of compact Stein surfaces}, preprint
(arXiv:math.GT/0103106).

\bibitem{DG}
F. Ding and H. Geiges
{\em Symplectic fillability of tight contact structures on torus bundles},
Alg. and Geom. Topology, {\bf1} (2001), 153--172.

\bibitem{DK} 
S. Donaldson and P. Kronheimer,
{\em Geometry of 4-manifolds},
Oxford Univ. Press, 1990.

\bibitem{E1}
Ya. Eliashberg,
{\em Filling by holomorphic disks and its applications},
Geometry of Low-Dimensional Manifolds: 2, Proc. Durham Symp. 1989, 
London Math. Soc. Lecture Notes, {\bf151}, Cambridge Univ. Press, 1990, 
45--67. 
 
\bibitem{E2}
Ya. Eliashberg,
{\em Unique holomorphically fillable contact structure on the 
3-torus},
Internat. Math. Res. Notices (1996), 77--82.

\bibitem{EH}
J. Etnyre and K. Honda,
{\em On the non-existence of tight structures},
preprint (arXiv:math.GT/9910115).

\bibitem{FS}
R. Fintushel and R. Stern,
{\em Knots, links and 4-manifolds},
Invent. Math. {\bf134} (1998), 363--400.


\bibitem{Fro} 
K. Fr\o yshov, {\em The Seiberg-Witten equations and
four-manifolds with boundary}, Math. Res. Lett. {\bf 3} (1996),
373-390.


\bibitem{F}
M. Furuta,
{\em Monopole equation and the $\frac{11}{8}$-conjecture},
Math. Res. Lett., to appear.

\bibitem{Gi}
E. Giroux, 
{\em Une structure de contact, meme tendue est plus ou moins tordue}, 
Ann. Scient. Ecole Normale Sup. {\bf27} (1994), 697--705. 

\bibitem{G}
R. Gompf,
{\em Handlebody construction of Stein surfaces}, Ann. of Math.
{\bf148} (1998), 619--693.


\bibitem{GS}
R. Gompf and A. Stipsicz,
{\em 4-manifolds and Kirby calculus},
Grad. Stud. in Math. AMS, 1999.

\bibitem{Gr}
H. Grauert
{\em On Levi's problem}, 
Ann. Math. {\bf68} (1958), 460--472. 

\bibitem{Ho}
K. Honda,
{\em personal communication}.

\bibitem{Ka}
Y. Kanda,
{\em The classification of tight contact structures on the 3-torus},
Comm. Anal. and Geom., {\bf5} (1997), 413--438.

\bibitem{L1}
P. Lisca,
{\em Symplectic fillings and positive scalar curvature}, 
Geom. Topol. {\bf2} (1998), 103--116. 

\bibitem{L2}
P. Lisca,
{\em On symplectic fillings of 3-manifolds}, 
Proc. of G\"okova Geometry and Topology Conference (1998), 151--160.

\bibitem{LM}
P. Lisca and G. Mati\'c,
{\em Tight contact structures and Seiberg-Witten invariants}, 
Invent. Math. {\bf129} (1997), 509--525.


\bibitem{Mc}
D. McDuff,
{\em The structure of rational and ruled symplectic 4-manifolds}
Journal  of the Amer. Math. Soc. {\bf3} (1990), 679--712.

\bibitem{Mo}
J. Morgan, 
{\em The Seiberg-Witten equations and applications to the topology
of smooth four-manifolds}, Math. Notes {\bf44}, Princeton University Press,
Princeton NJ, 1996.

\bibitem{MSz}
J. Morgan and Z. Szab\'o,
{\em Homotopy K3-surfaces and mod 2 Seiberg-Witten invariants},
Math. Res. Lett. {\bf4} (1997), 17--21.

\bibitem{MOY}
T. Mrowka, P. Ozsv\'ath and B. Yu,
{\em Seiberg-Witten monopoles on Seifert fibered spaces},
Comm. Anal. Geom. {\bf5} (1997), 685--791.

\bibitem{OO}
H. Ohta and K. Ono,
{\em Simple singularities and topology of symplectically filling
4-manifold},
Comment. Math. Helv. {\bf74} (1999), 575--590.

\bibitem{OSz}
P. Ozsv\'ath and Z. Szab\'o,
{\em The symplectic Thom conjecture}, Annals of Math.
{\bf152} (2000), 93--124.


\bibitem{OSz2}
P. Ozsv\'ath and Z. Szab\'o,
{\em On embedding of circle bundles in 4-manifolds},
Math. Res. Letters {\bf7} (2000), 657--669.

\bibitem{Ro} 
V. Rohlin, 
{\em A three-dimensional manifold is the boundary of a four-dimensional one}, 
Dokl. Akad. Nauk SSSR {\bf81} (1951), 355--357.

\bibitem{Sgok}
A. Stipsicz,
{\em Gauge theory and Stein fillings of certain 3-manifolds},
in preparation.

\bibitem{SSz}
A. Stipsicz and Z. Szab\'o,
{\em Gluing 4-manifolds along $\Sigma (2,3,11)$},
Top. and its Applications, {\bf106} (2000), 293--304.

\bibitem{T}
C. Taubes,
{\em The Seiberg-Witten invariants and 4-manifolds with essential 
tori}, Geometry and Topology {\bf5} (2001), 441--519.


\bibitem{W}
E. Witten,
{\em Monopoles and four-manifolds}, Math. Res. Lett. {\bf1} (1994), 769--796.


\end{thebibliography}
\end{document}